\documentclass[11pt]{article}
\usepackage{fancyhdr}
\usepackage{isomath}
\usepackage{amsmath}
\usepackage{amsbsy}
\usepackage{amssymb}
\usepackage{amscd}
\usepackage{amsfonts}
\usepackage{graphicx,color}
\usepackage{verbatim}
\usepackage{euscript}
\usepackage{alltt}
\usepackage{stmaryrd}
\usepackage{subfigure}
\usepackage{relsize}
\usepackage{enumitem}
\usepackage{url}
\usepackage[stable]{footmisc}

\DeclareGraphicsExtensions{.eps,.pdf}

\usepackage[margin=1in]{geometry}

\newcommand{\p}{\partial}
\newcommand{\dee}{\mathcal{D}}
\newcommand{\scl}{\mathcal{L}}

\date{}
\begin{document}
\title{A dual variational principle for nonlinear dislocation dynamics}

\author{Amit Acharya\thanks{Department of Civil \& Environmental Engineering, and Center for Nonlinear Analysis, Carnegie Mellon University, Pittsburgh, PA 15213, email: acharyaamit@cmu.edu.}}
\maketitle
$\qquad \quad$ \textit{This paper is dedicated to Roger Fosdick on the occasion of his 85th birthday}
\begin{abstract}
\noindent A dual variational principle is defined for the nonlinear system of PDE describing the dynamics of dislocations in elastic solids. The dual variational principle accounting for a specified set of initial and boundary conditions for a general class of PDE is also developed.
\noindent 
\end{abstract}

\section{Introduction}
	 Unlike the physics of the microscopic structure of sub-atomic particles (e.g. `core' of an electron), much more is physically known about the microscopic structure of dislocations and their mutual interactions, as well as their interactions with applied loads, within a (nonlinear) elastic crystal, both through direct experimental observation and through lattice statics\-/\-molecular dynamics/density functional theory calculations.
	 %%Unlike the physics of the microscopic structure of sub-atomic particles (e.g. `core' of an electron), much is physically known, through direct experimental observation and lattice statics/molecular dynamics/density functional theory calculations about the microscopic structure of dislocations and their mutual interactions, as well as with applied loads through boundary conditions, within a (nonlinear) elastic crystal. 
	 Due to this knowledge, physically well-justified and transparent mathematical models can be posited for the phenomena, with the possibility of systematic refinement to include more detail when deemed necessary after mathematical study and comparison with experiment. There is a long and distinguished history of the study of dislocations in elasticity in the classical setting, see, e.g., \cite{weingarten, volterra, peierls, nabarro, hirth_lothe}, the continuously distributed setting, e.g., \cite{nye, kroner, mura}, \cite[including second-order effects]{willis} and \cite{fox}, and the connections of some of the kinematic aspects of dislocations to non-Riemannian Geometry \cite{kondo, kroner, bilby}. As well, techniques for developing well-set, classical thermomechanical theories of the mechanics of continuous media comprising different types of materials exhibiting strongly nonlinear behavior and satisfying the relevant invariances and material symmetries are available \cite{cft,nlftm, coleman_noll, coleman_gurtin} and \cite{fox}. These ideas and techniques have been synthesized and extended to produce the theory/model of dislocation mechanics stated in \cite{fdm}, as reviewed in \cite{action}. The theory admits the minimal specification of an energy density function $\psi(W)$, where $W$ is the inverse elastic distortion field (not necessarily a gradient), and that of a dislocation velocity field, the function $V_s(\alpha, W, \rho)$ in \eqref{eq:FDM}, which, when guided by the requirements of being proportional to its derived thermodynamic driving force, is a specified function of the thermodynamically derived Cauchy stress tensor $T_{ij} = - \rho W_{ki} \psi'_{kj}$ and the dislocation density tensor $\alpha_{ij}$, admitting a scalar or matrix of material constants representing dislocation mobility. Here, $\psi'_{ij} = \p_{W_{ij}} \psi$, and it suffices to use a rectangular Cartesian coordinate system and tensor components w.r.t its basis in this Section. The time variable is represented by the symbol $t$ and not used as an index.

For prescribed static dislocation fields the framework is shown to be able to compute the stress and energy fields of such distributions in bodies of arbitrary geometry and general elastic symmetries \cite{arora_2, arora_1}. Similarly for prescribed dislocation velocity field, the setup is shown to be able to compute the evolution of the dislocation field \cite{arora_1}. And the evolution in the fully coupled case also has been shown to work well to predict nonsingular dislocation cores, dislocation annihilation, dissociation and stress-mediated interaction when restricted to dislocation motion within a planar layer in a 3--d body \cite{zhang2015single} within a `small deformation' ansatz.

The phenomenon of macroscopic plasticity of crystalline materials corresponds to the collective dynamical behavior of a very large number of dislocation curves in an elastic body under generally time-dependent loads. While experimental observations and real practical applications of plasticity abound, it is fair to say that there does not exist a fundamental theory that arises as a coarse-graining of nonlinear dislocation dynamics as described above (or by any other model). The phenomenon of plasticity shows fascinating dynamical changes as a function of initial conditions and tamely evolving driving loads - e.g., yielding, Stage I, II, III, IV behaviors as a function of applied load temperature and initial crystal orientation, intricate patterned dislocation microstructure formation such as cells and sub-grain boundaries to name only a few - with no established fundamental theory for understanding them (the phenomenon is even richer, with rapidly driven situations also being of theoretical and practical interest). It is in this context that we would like to use a path integral implementation of the dynamics represented by \eqref{eq:FDM} to evaluate how much of the reality of macroscopic plasticity can be understood by the combination of the model and the technique. The rough expectation is to be able to interpret drastic changes of overall behavior observed in reality as statistical phase transitions as understood in Effective/Statistical Field Theory.

A first step in this program is to define an action functional for the system \eqref{eq:FDM} which, in the first instance, does not emanate from a variational principle; it is this objective that is tackled in Secs. \ref{sec:fin_dim}, \ref{sec:dual_finite_disloc}, and \ref{sec:dual_lin_disloc}, refining the work in \cite{action} following the ideas in \cite{action_2}. Section \ref{sec:nonlin_bc} develops the variational principle accounting for a specified set of initial and boundary conditions for a nonlinear system of second-order PDE expressed as a first-order system.

%%{\color{red}%% 
Variational principles for `small deformation,' static dislocation mechanics and internal stress problems by the method of `eigenstrains' is presented in \cite{kroner_2, berdich}, and there is a modern literature involving rigorous analysis reviewed, in detail in \cite{arora_2}. Our work involves finite-deformation, nonlinear dislocation dynamics including inertia.%%}

	\section{The essential idea: An optimization problem for an algebraic system of equations\footnote{I thank Vladimir Sverak for insisting on the `simplest,' transparent explanation of the ideas in \cite{action, action_2}. This brief Section is a result of that effort.}}\label{sec:fin_dim}
	Consider a generally nonlinear system of algebraic equations in the variables $x \in \mathbb{R}^n$ given by
	\begin{equation}\label{eq:alg_sys}
	    A_\alpha (x) = 0,
	\end{equation}
	where $A: \mathbb{R}^n \to \mathbb{R}^N$ is a given function (a simple example would be $A_\alpha (x) = \bar{A}_{\alpha i} \, x^i - b_\alpha$, $\alpha = 1 \ldots N, i = 1 \ldots n$, where $\bar{A}$ is a constant matrix, \textit{not necessarily symmetric} (when $n = N$), and $b$ is a constant vector). We allow for all possibilities $0 < n \lesseqqgtr N > 0$.
	
	The goal is to construct an objective function whose critical points solve the system \eqref{eq:alg_sys} (when a solution exists) by defining an appropriate $x^* \in \mathbb{R}^n$ satisfying  $A_\alpha (x^*) = 0$.
	
	For this, consider first the auxiliary function
	\begin{equation*}
	    \widehat{S}_H(x,z) = z^\alpha A_\alpha (x) + H(x)
	\end{equation*}
	(where $H$ belongs to a class of scalar-valued function to be defined shortly) and define
	\begin{equation*}
	    S_H(z) = z^\alpha A_\alpha(x_H (z)) + H(x_H(z))
	\end{equation*}
	with the requirement that the system of equations
	\begin{equation}\label{eq:H_fin_dim}
	    z^\alpha \frac{\p A_\alpha}{\p x^i}(x) + \frac{\p H}{\p x^i}(x) = 0
	\end{equation}
	be solvable for the function $x = x_H(z)$ through the choice of $H$, and \textit{any} function $H$ that facilitates such a solution qualifies for the proposed scheme. 
	
	In other words, given a specific $H$, it should be possible to define a function $x_H(z)$ that satisfies 
	\begin{equation*} 
	z^\alpha \p_{x^i} A_\alpha (x_H(z)) + \p_{x^i} H (x_H(z)) = 0 \quad \forall z \in \mathbb{R}^N
	\end{equation*}
	(the domain of the function $x_H$ may accommodate more intricacies, but for now we stick to the simplest possibility). Note that \eqref{eq:H_fin_dim} is a set of $n$ equations in $n$ unknowns regardless of $N$ ($z$ for this argument is a parameter).
	
	Assuming this is possible, we have
	\begin{equation*}
	    \frac{\p S_H}{\p z^\beta} (z) = A_\beta(x_H(z)) +  \left( z^\alpha \frac{\p A_\alpha}{\p x^i}(x_H(z)) + \frac{\p H}{\p x^i}(x_H(z)) \right) \frac{\p x^i_H}{\p z^\beta}(z) = A_\beta(x_H(z)),
	\end{equation*}
	using \eqref{eq:H_fin_dim}. Thus,
	\begin{itemize}[noitemsep,topsep=2pt]
	    \item if $z_0$ is a critical point of the objective function $S_H$ satisfying $\p_{z^\beta} S_H(z_0) = 0$, then the system $A_\alpha(x) = 0$ has a solution defined by $x_H(z_0)$; 
	    \item if the system $A_\alpha(x) = 0$ has a unique solution, say $y$, and if $z^H_0$ is any critical point of $S_H$, then $x_H\left(z^H_0 \right) = y$, for all admissible $H$.
	    \item If $A_\alpha(x) = 0$ has non-unique solutions, but $\p_{z^\beta} S(z) = 0$ ($N$ equations in $N$ unknowns) has a unique solution for a specific choice of the function $z \mapsto x_H(z)$ related to a choice of $H$, then such a choice of $H$ may be considered a selection criterion for imparting uniqueness to the problem $A_\alpha(x) = 0$.
        \item Finally, to see the difference of this approach with the Least-Squares (LS) Method, for a linear system $\bar{A} x = b$, the LS governing equations are given by
        \[
        \bar{A}^T \bar{A} z = \bar{A}^T b,
        \]
        with LS solution defined as $z$ even when the original problem $\bar{A} x = b$ does not have a solution (i.e., when $b$ is not in the column space of $\bar{A}$). The LS problem always has a solution, of course. In contrast, in the present duality-based approach with quadratic $H(x) = \frac{1}{2} x^T x$ the governing equation is
        \[
        \bar{A}\bar{A}^T z = b
        \]
        with solution to $\bar{A} x = b$ given by $x = \bar{A}^T z$, and the problem has a solution only when $\bar{A} x = b$ has a solution, since the column spaces of the matrices $\bar{A}$ and $\bar{A}\bar{A}^T$ are identical.

        As a practical matter, the latter approach appears to have, at least in principle, advantages for solving large, \textit{consistent}, underdetermined systems as the size of the matrix $\bar{A}\bar{A}^T$ is much smaller than that of $\bar{A}^T \bar{A}$ in this situation, with due consideration paid to conditioning-related robustness issues (cf. \cite{george1984solution, heath2018scientific}, \cite[pp. 299-300]{golub2013matrix}).
	\end{itemize}
	\section{A class of variational principles for nonlinear dislocation mechanics}\label{sec:dual_finite_disloc}
	We implement the idea of Sec.~\ref{sec:fin_dim} to define an action(s) for the nonlinear partial differential equations of dislocation mechanics given by
	\begin{equation}\label{eq:FDM}
	\begin{aligned}
	    0 & = e_{jrs} \p_r W_{is} + \alpha_{ij} \\
	    0 & = \p_t W_{ij} + \p_j (W_{ik} v_k) - v_k e_{rkj} \alpha_{ir} - e_{jrs} \alpha_{ir} V_s(\alpha, W, \rho) = \p_t W_{ij} + v_k \p_k W_{ij} + W_{ik} \p_j v_k - e_{jrs} \alpha_{ir} V_s(\alpha,W, \rho)\\
	   0 & =  \p_t \rho + \p_k (\rho v_k)\\
	   0 & = \p_t (\rho v_i) + \p_j (\rho v_i v_j) + \p_j (\rho W_{ki} \psi'_{kj}).
	\end{aligned}
	\end{equation}
	The physical basis of this system of equations is explained in \cite{action}. Briefly, the first equation is the equation of elastic incompatibility. The second reflects the compatibility between the velocity gradient, the rate of the change of the elastic distortion and the rate of permanent deformation produced by the motion of dislocations. The third equation is balance of mass, and the fourth, the balance of linear momentum. Setting $\alpha = 0$ in the system above gives the equations of nonlinear elasticity written in an Eulerian setting.
	
	First define the functional
	\begin{equation*}
	\begin{aligned}
	    \widehat{S}_H[A,W,\theta,\rho,\lambda,v,B,\alpha] 
	     = \int_{[0,T]\times \Omega} \, dt d^3x & \, - W_{ij} \p_t A_{ij} - W_{ik} v_k \p_j A_{ij} - A_{ij} v_k e_{rkj} \alpha_{ir} - A_{ij} e_{jrs} \alpha_{ir} V_s(\alpha, W, \rho) \\
	    & - \rho \p_t \theta - \rho v_k \p_k \theta \\
	    & - \rho v_i \p_t \lambda_i - \rho v_i v_j \p_j \lambda_i - \rho W_{ki} \psi'_{kj} \p_j \lambda_i\\
	    & - e_{jrs} W_{is} \p_r B_{ij} + B_{ij} \alpha_{ij}\\
	    & + H(W, \rho, v, \alpha),
	\end{aligned}
	\end{equation*}
	which is obtained by converting \eqref{eq:FDM} to scalar form by taking inner products with the `dual' fields 
	\[
	D = (A, \theta, \lambda, B),
	\]
	integrating by parts on the space-time domain assuming the dual fields vanish on the boundary of the domain, and adding the potential $H$.
	Now define
	\begin{equation*}
	    U := (W, \rho, v, \alpha) \qquad \mbox{and} \qquad \dee := (\p_t A, \nabla A, A, \p_t \theta, \nabla \theta, \p_t \lambda, \nabla \lambda, \nabla B, B)
	\end{equation*}
	(note `$\dee \neq D$') and require that there exists a function
	\begin{equation}\label{eq:Uofdee}
	    U_H(\dee) = (W_H(\dee),\rho_H(\dee),v_H(\dee), \alpha_H(\dee))
	\end{equation}
	such that \textit{for the functional $S_H[A, \theta, \lambda, B]$ of the dual fields defined as}
	\begin{equation}\label{eq:SH}
	   \int_{[0,T]\times \Omega} \, dt d^3x \ \scl_H(\dee, U_H(\dee)) = S_H [A, \theta, \lambda, B] := \widehat{S}_H [A, W_H(\dee), \theta, \rho_H(\dee), \lambda, v_H(\dee), B, \alpha_H(\dee)],
	\end{equation}
	the first variation is given by (we suppress the subscript $_H$ on the elements of $U_H$ for notational simplicity)
	\begin{equation}\label{eq:var_SH}
	\begin{aligned}
	     \delta S_H = \int_{[0,T]\times \Omega} \, dt d^3x & \, - W_{ij}(\dee) \p_t \delta A_{ij} - W_{ik} (\dee) v_k(\dee) \p_j \delta A_{ij} - \delta A_{ij} v_k (\dee) r_{rkj} \alpha_{ir} (\dee) \\
	     & \qquad \qquad \qquad \qquad - \delta A_{ij} e_{jrs} \alpha_{ir} (\dee) V_s (\alpha(\dee), W(\dee), \rho(\dee))\\
	     & - \rho(\dee) \p_t \delta \theta - \rho(\dee) v_k (\dee) \p_k \delta \theta \\
	     & - \rho(\dee) v_i(\dee) \p_t \delta \lambda_i - \rho(\dee) v_i(\dee) v_j(\dee) \p_j \delta \lambda_i - \rho(\dee) W_{ki}(\dee) \psi'_{kj}(W(\dee)) \p_j \delta \lambda_i\\
	     & - e_{jrs} W_{is} (\dee) \p_r \delta B_{ij} + \delta B_{ij} \alpha_{ij}(\dee),
	\end{aligned}
	\end{equation}
	\textit{a condition that is satisfied if the system}
	\begin{equation}
	    \begin{aligned}\label{eq:invert}
	        \frac{\p \scl_H}{\p W_{lp}} = - \p_t A_{lp} - v_p \p_j A_{lj} - A_{ij} e_{jrs} \alpha_{ir} \frac{\p V_s}{\p W_{lp}} (\alpha, W, \rho) - e_{jrp} \p_r B_{lj} & \\
	        - \rho \left( \psi'_{lj}(W) \p_j \lambda_p \right. \left. + \ W_{ki} \psi''_{kjlp} \p_j \lambda_i \right)& + \frac{\p H}{\p W_{lp}}(W,\rho, v, \alpha) = 0\\
	       \frac{\p \scl_H}{\p \rho} = - A_{ij} e_{jrs} \alpha_{ir} \frac{\p V_s}{\p \rho} (\alpha, W, \rho) - \p_t \theta - v_k \p_k \theta - v_i \p_t \lambda_i - v_i v_j \p_j \lambda_i & - W_{ki} \psi'_{kj} \p_j \lambda_i\\
	       & + \frac{\p H}{\p \rho} (W, \rho, v, \alpha) = 0\\
	        \frac{\p \scl_H}{\p v_p} = - W_{ip} \p_j A_{ij} - A_{ij} e_{rpj} \alpha_{ir} - \rho \p_p \theta - \rho \p_t \lambda_p - \rho v_j \p_j \lambda_p & - \rho v_i \p_p \lambda_i \\
	        & + \frac{\p H}{\p v_p} (W, \rho, v, \alpha) = 0\\
	        \frac{\p \scl_H}{\p \alpha_{lp}} = - A_{lj} v_k e_{pkj} - A_{lj} e_{jps} V_s(\alpha, W, \rho) - A_{ij} e_{jrs} \alpha_{ir} \frac{\p V_s}{\p \alpha_{lp}} (\alpha, W, \rho) & + B_{lp} \\
	        & + \frac{\p H}{\p \alpha_{lp}}(W,\rho, v, \alpha) = 0,
	    \end{aligned}
	\end{equation}
	\textit{can be solved in the form of}
	\begin{equation*}
	   (W,\rho, v, \alpha) = U_H(\dee).
	\end{equation*}
	This is so, since solving \eqref{eq:invert} defines $U_H(\dee)$ that ensures $\frac{\p \scl_H}{\p U} (\dee, U_H(\dee)) = 0$ which then implies
	\begin{equation*}
	    \frac{\p \scl_H}{\p U} (\dee, U_H(\dee)) \cdot \frac{\p U_H}{\p \dee} (\dee) \cdot \delta \dee = 0 \qquad \mbox{for all} \ \dee.
	\end{equation*}
     Note that \eqref{eq:var_SH} then is simply
	\begin{equation*}
	    \delta S_H =  \int_{[0,T]\times \Omega} \, dt d^3x\  \dfrac{\p \scl_H}{\p \dee} \cdot \delta \dee,
	\end{equation*}
	and requiring 
	\begin{equation*}
	    \delta S_H = 0 \ \mbox{for all variations}\  \delta D \ \mbox{that vanish on the boundary of} \  \Omega \times [0,T]
	\end{equation*} 
	  shows that the Euler-Lagrange equations of the functional $S_H$ defined in \eqref{eq:SH} are the equations of \eqref{eq:FDM} with the substitution
	\begin{equation*}
	    (W,\rho,v, \alpha) = U_H(\dee).
	\end{equation*}
	This is so because $\scl_H$ is necessarily linear in its first argument, see \eqref{eq:SH}.
	
	To summarize, \textit{the primal equations \eqref{eq:FDM} of dislocation mechanics are the Euler-Lagrange equations of any of the dual functionals, written in terms of particular specific combinations (mappings) of the dual fields for each choice of the function $H$, each specific mapping defining the primal fields. Thus one may think of the primal fields as “gauge invariant” observable combinations of the dual fields (“gauge fields”) satisfying one specific set of equations (the primal system). While this is not how gauge fields appear in traditional gauge theories of physics, it is interesting that a completely different starting point and approach raise somewhat similar invariance structures that may be interpreted as symmetries}.
	
	As for the plausibility of being able to solve the \textit{algebraic} system \eqref{eq:invert} given a specific $\dee$, consider $H$ to be separately quadratic in each of its arguments, say $U_A$, with large in magnitude coefficient, so, e.g., $H = \frac{1}{2} \alpha_W W_{ij} W_{ij} + \cdots$, with $ 1 \ll |\alpha_W|$. Then, assuming the solution of the Euler-Lagrange equations are bounded in some appropriate sense, \eqref{eq:invert} can indeed be solved to define $U_H(\dee)$, and it has to be made sure that the solutions of the Euler-Lagrange equations (using this function) indeed satisfy the assumed bounds. To ensure that this latter condition is satisfied one has a large class of $H$ functions to operate with but, at any rate, this is a delicate question of analysis, including how the requirement can be relaxed if required, in the context of solving the dual variational problem (and not necessarily its Euler-Lagrange equations).
	
	We end this section with the following remarks:
	\begin{itemize}[noitemsep, topsep = 2 pt]
	    \item Our system \eqref{eq:FDM} does not involve multi-valued fields or non-simply connected domains for defining dislocation dynamics, but is fully capable of representing the topological charge of dislocation lines with its ingredients.
	    \item Based on the explorations of stress-coupled dislocation motion presented in \cite{ach_tartar, zhang2015single}, the `primal' system requires a `core-energy' in the form of the dependence of the energy function $\psi$ on the dislocation density $\alpha$ as well. This results in the dislocation velocity depending on the $\mbox{curl}\, \alpha$. Such a dependence is accommodated within our `action-generating' scheme by adding an extra variable and equation to the system \eqref{eq:FDM} of the form $e_{jrs} \p_r \alpha_{is} = \beta_{ij}$ and writing the dislocation velocity as $V_s = V_s(\alpha, W, \rho, \beta)$. This would have the effect of increasing the number of fields in the dual problem as well.
	    
	    It is an interesting question whether the precise definition of a formally `small' core energy contribution with a small parameter representing microscopic physics can make a difference in the development of an accurate model for the prediction of macroscopic behavior, and whether such a device should be allowed in the class of models admitted. Physically, in the context of the physics of dislocation dynamics, there appears to be no reason to exclude the possibility of the importance of such effects and, in fact, allows more precise physics to be incorporated in the description of gross macroscopic behavior (which is, admittedly, a double-edged sword in the context of coarse-graining). Some evidence to support such an expectation is also provided by the mathematically rigorous study of the inviscid Burgers equation, `regularized' by a small viscous effect in one case and by dispersion in another \cite{lax1973hyperbolic, lax1986dispersive, lax1979zero}.
	    
	    Based on the above observation, one advantage of the `dual' formalism proposed herein may be that when the microscopic physics to be added is not even qualitatively understood with certainty, working with a regularization on the dual side, may be guided solely by the aim of producing a `good' dual extremal, i.e., with guaranteed existence in an appropriate function space. Doing so appears to require no modification to the physics of the primal problem, and then  the limit of dual solutions, as the regularization parameter vanishes, may be studied. 
	    \item In the context of an action functional that simply has as its Euler-Lagrange equation the given system of PDEs, the proposed scheme delivers, at least formally and under the stated requirements, what is needed. However, if the action functional is to be used in a path integral, dual fields $D$ other than extremals matter as well. In this sense it is reasonable to demand that the added potential $H$ in $\scl_H$ be subject to further requirements of invariance that may obstruct the inversion process required to define the function $U_H(\dee)$. In case such a restriction is so severe as not to allow the definition of even a single `change of variables' ($U_H(\dee)$), through the choice of some $H$), one can retain both the fields $W,A$ and still obtain a relevant action functional, as shown in \cite{action}.
	\end{itemize}
\section{Linear dislocation mechanics}\label{sec:dual_lin_disloc}
We illustrate the proposed technique with a very closely related one (using a Legendre transform, cf.~\cite{zaanen2021crystal, action_2}) in the simplified setting of linear dislocation mechanics with a prescribed dislocation velocity field $V$ in space-time along with the ansatz
\begin{equation*}
    \begin{aligned}
    U_{ij} & := \delta_{ij} - W_{ij} \\
    T_{ij} & := C_{ijkl} U_{kl},
    \end{aligned}
\end{equation*}
ignoring all nonlinearities in \eqref{eq:FDM} and assuming the mass density field $\rho$ to a be specified field. The ansatz is justified for small elastic distortions $(U)$ about the ground state (cf., \cite{action}). We note that $C_{ijkl}$ is necessarily symmetric in $(k,l)$ and $(i,j)$ so that it is not invertible on the space of all second order tensors (and hence the stress only depends on the elastic strain, the symmmetric part of $U$). With these assumptions, the system \eqref{eq:FDM} may be expressed as
\begin{equation}\label{eq:lin_FDM}
    \begin{aligned}
    0 & = \p_t (\rho v_i) - \p_j (C_{ijkl} U_{kl})\\
    0 & = \p_j v_i - \p_t U_{ij} - e_{jrs} \alpha_{ir} V_s \\
    0 & = e_{jrs} \p_r U_{is} - \alpha_{ij}.
    \end{aligned}
\end{equation}
Taking inner products of these equations with the dual fields $D = (\lambda, A, B)$ that vanish on the boundary and utilizing an arbitrary function $M$ convex in the list of arguments $M(v, U, \alpha)$ we define the functional 
\begin{equation*}
\begin{aligned}
   \widehat{S}[A,U,B, \alpha, \lambda, v] = \int_{[0,T]\times \Omega} \, dt d^3x & \ v_i (- \p_j A_{ij} - \rho \p_t \lambda_i) \\
   & \ + U_{ij} ( \p_t A_{ij} - e_{sjr} \p_r B_{is} + C_{ijkl} \p_l \lambda_k) \\
   & \ + \alpha_{ir} (- A_{ij} e_{jrs} V_s - B_{ir}) - M(U, \alpha, v)
\end{aligned}
\end{equation*}
Defining
\begin{equation*}
    p := (-\p_j A_{ij} - \rho \p_t \lambda_i, \p_t A_{ij} - e_{sjr} \p_r B_{is} + C_{ijkl} \p_l \lambda_k, - A_{ij} e_{jrs} V_s - B_{ir}); \qquad \qquad Q := (v, U, \alpha)
\end{equation*}
and $M^*(p)$ the Legendre transform of $M(Q)$ given by
\begin{equation}
    \begin{aligned}\label{eq:lin_map}
    (v_M(p), U_M(p), \alpha_M(p)) =: Q_M(p) &= (\p_Q M)^{-1}(p)\\
    M^*(p) & = Q_M(p) \cdot p - M( Q_M(p))\\
    \p_p M^*(p) & = Q_M(p)
    \end{aligned}
\end{equation}
(well-defined because of the convexity of $M(Q)$), we define the dual action, $S_M[D]$,
\begin{equation*}
    \widehat{S}[A,U_M(p), B, \alpha_M(p), \lambda , v_M(p)] =: S_M[D] = \int_{[0,T]\times \Omega} \, dt d^3x  \ M^*(p)
\end{equation*}
whose first variation is given by (after an integration by parts)
\begin{equation*}
    \begin{aligned}
    \delta S_M & = \int_{[0,T]\times \Omega} dt d^3x \  Q_M(p)\, \delta p \\
    & = \int_{[0,T]\times \Omega}  dt d^3x \  \delta \lambda_i \left( \p_t(\rho v_i(p)) - \p_j (C_{ijkl} U_{kl}(p) ) \right) \\
    &  \qquad  \qquad \qquad \, + \ \delta A_{ij} \left(\p_j (v_i(p)) - \p_t (U_{ij}(p)) - e_{jrs} \alpha_{ir}(p) V_s \right)  \\
    &  \qquad \qquad  \qquad \, + \ \delta B_{is} ( e_{srj} \p_r (U_{ij}(p)) - \alpha_{is}(p) ),
    \end{aligned}
\end{equation*}
(where we have dropped the subscript $_M$ on the dual-to-primal mapping fields for notational convenience). Thus, the dual Euler-Lagrange equations are the system \eqref{eq:lin_FDM} expressed in terms of the dual fields through the mapping codified in \eqref{eq:lin_map}$_2$, \textit{regardless of the convex potential} $M$ \textit{chosen to define the dual functional} $S_M$.

This exercise exposes an interesting fact in a simple setting. Clearly, for $M$ to be convex in $U$ it cannot be invariant as it has to depend on the skew-symmetric part of the latter - and rotational invariance/invariance under superposed rigid deformations in the linear setting precludes such a dependence. However, the use of such a potential in the dual theory does not in any way obstruct the definition of correct physics as embodied in the Euler-Lagrange equations solved.

\section{Dual variational principle for a primal problem with initial and boundary conditions}\label{sec:nonlin_bc}
Consider the system of PDE:
\begin{equation}\label{eq:H-J_ibc}
    \begin{aligned}
    \p_t u_I &= \mathbb{A}_{IJ}\, u_J + \mathbb{B}_{IJk}\, \p_k u_{J} + \mathbb{C}_{IJkl} \, \p_l B_{Jk} + \mathfrak{f}_I(u, B, C) + \p_k ( \mathfrak{A}_{Ik}(u, B, C))\\
    \p_i u_I & = B_{Ii} \\
    \p_j B_{Ii} & = C_{Iij},
    \end{aligned}
\end{equation}
where $\mathbb{A}, \mathbb{B}, \mathbb{C}$ are arrays of real constants, $\mathfrak{f}_I$ and $\mathfrak{A}_{Ij}$ are, for each $I,j$, given, real-valued, smooth functions of the arguments shown, uppercase Latin indices span $1$ to $n$, and lowercase Latin indices representing space-dimensions span $1$ to $1 \leq d \leq 3$, and $t$ is time.  

The functions $\mathfrak{f}, \mathfrak{A}$ do not contain any terms linear in the array $(u, B, C)$.

Let the initial and boundary conditions for \eqref{eq:H-J_ibc} be
\begin{equation}\label{eq:ibc}
    \begin{aligned}
    & u_I(x, 0)  = \overline{u}^{(i)}_I(x),  \qquad  x \in \Omega \\
    & \mathfrak{A}_{Ik}(u(x,t),B(x,t),C(x,t)) \, n_k(x)  \\
    & \qquad + \left( \mathbb{B}_{IJk} u_J(x,t) + \mathbb{C}_{IJik} B_{Ji}(x,t)  \right) n_k(x) = \overline{\tau}_I(x,t),  \qquad   (x,t) \in \p \Omega_{\tau}(x) \times (0,T]\\
   &  u_I(x,t)  = \overline{u}^{(b)}_I (x,t), \qquad (x,t) \in \p \Omega_{u}(x) \times (0,T]\\
    & \p_i u_I (x, t)  = B_{Ii}(x,t) = \overline{B}_{Ii} (x,t), \qquad  (x,t) \in \p \Omega_{\nabla u}(x) \times (0,T]
    \end{aligned}
\end{equation}
where $n$ represent the outward unit normal to the boundary of the domain $\p \Omega$, the functions with overhead bars are prescribed, and the subsets $\p \Omega_{\tau}, \p \Omega_{u}, \p \Omega_{\nabla u}$ of the spatial boundary $\p \Omega$ can very well be empty for a specific problem.

The initial and boundary conditions \eqref{eq:ibc} for the primal system \eqref{eq:H-J_ibc} is simply a set of conditions that encompass the commonly encountered ones for up to second-order systems of partial differential equations; the present work does not deal with the question of well-posedness of the system \eqref{eq:H-J_ibc} with this specified set of initial and boundary conditions. For instance, it may very well be that in a specific problem, well-possedness requires  $\p_i u_I = B_{Ii}$ not to be specified on any part of the boundary of $\Omega$. In that case $\p\Omega_{\nabla} u$ can be chosen to be the empty set, and then, as suggested by \eqref{eq:S_ibvp} below, the dual field $\rho$ needs to be prescribed on the entirety of $\p \Omega$ for all times.

Our scheme \cite{action_2} then suggests defining
\begin{equation*}
\begin{aligned}
    \widehat{S}[u,B,C,\lambda,\gamma,\rho]  = \int_{\Omega \times (0,T)} d^3x dt  \, & - \p_t \lambda_I u_I - \lambda_I \mathfrak{f}_I(u, B, C) + \mathfrak{A}_{Ik}(u, B, C) \p_k \lambda_I \\
   &  - u_I \mathbb{A}_{JI} \lambda_J + u_I \, \mathbb{B}_{JIk} \, \p_k \lambda_J + B_{Ii}\, \mathbb{C}_{JIil}\, \p_l \lambda_J\\
     & - u_I \p_i \gamma_{Ii} - \gamma_{Ii} B_{Ii} \\
     &  - B_{Ii} \p_j \rho_{Iij} - \rho_{Iij} C_{Iij} + H(u, B, C, x, t),
\end{aligned}
\end{equation*}
\emph{where we have allowed for $H$ to depend on $(x,t)$ as well. This generality is practically useful, especially when $H$ plays the role of a selection criterion for non-unique solutions of the primal problem}.

Next define the arrays $U, P, L, F$ and the function $M$:
\begin{equation*}
\begin{aligned}
    & U := (u_I, B_{Ii}, C_{Iij}), \qquad L := (\lambda_I, \ - \p_k \lambda_I) \\
    & P := (\p_t \lambda_I + \p_i \gamma_{Ii} + \mathbb{A}_{JI} \lambda_J - \mathbb{B}_{JIk} \p_k \lambda_J, \quad \gamma_{Ii} + \p_j \rho_{Iij} - \mathbb{C}_{JIil} \, \p_l \lambda_J, \quad \ \rho_{Iij}), \qquad F := (\mathfrak{f}_I, \mathfrak{A}_{Ij}) \\
    & M(U, L, x,t) := H(U, x, t) - L \cdot F(U)\\
\end{aligned}
\end{equation*}
(noting that $F$ is indeed a function of only $U$). 

Following \cite{action_2}, we now ask that $M$ (enabled by the choice of $H$) be such that 
\[ 
\exists \  U_H(P, L, x,t) \quad \mbox{satisfying} \quad \p_U M (U_H(P,L,x,t), L,x,t) = P, \ \forall P
\]
and in terms of this define
\begin{equation*}
    M^*(P,L,x,t) := U_H(P,L,x,t) \cdot P - M(U_H(P,L,x,t), L,x,t).
\end{equation*}
We then replace $U$ in $\widehat{S}$ by $U_H(P,L,x,t)$ to define the \textit{dual} functional
\begin{equation*}
\begin{aligned}
    S[\lambda, \gamma, \rho] & := \int_{\Omega \times (0,T)} d^3x dt  \  - P \cdot U_H(P,L,x,t) + M(U_H(P,L,x,t), L,x,t) \\
    & = - \int_{\Omega \times (0,T)}  d^3x dt \  M^*(P,L, x,t)
    \end{aligned}
\end{equation*}
whose first variation is given by
\begin{equation*}
    \delta S = \int_{\Omega \times (0,T)} d^3x dt  \, - \p_P M^* \cdot \delta P - \p_L M^* \cdot \delta L.
\end{equation*}
Now,  a calculation detailed in \cite[Sec.~6]{action_2}, shows that
\begin{equation*}
    \p_P M^*(P, L,x,t) = U_H(P,L,x,t), \qquad \p_L M^*(P,L,x,t) = F(U_H(P,L,x,t)).
\end{equation*}
Then, using the notation $U_H(P,L,x,t) = \left( u^H_I(P, L,x,t), \ B^H_{Ii}(P,L,x,t), \  C^H_{Iij}(P,L,x,t) \right), I = 1 \ \mbox{to} \ n;\  i,j = 1 \ \mbox{to} \ d$,
\begin{equation*}
    \begin{aligned}
    \delta S & = \int_{\Omega \times (0,T)} d^3x dt  \, - U_H(P,L,x,t) \cdot \delta P - F(U_H(P,L,x,t)) \cdot \delta L\\
    & = \int_{\Omega \times (0,T)} d^3x dt  \, - u^H_I(P,L,x,t) (\p_t \delta \lambda_I + \p_i \delta \gamma_{Ii} + \mathbb{A}_{JI} \delta \lambda_J - \mathbb{B}_{JIk} \p_k \delta \lambda_J ) \\
    & \qquad \qquad \qquad \qquad - B^H_{Ii}(P,L,x,t) (\delta \gamma_{Ii} + \p_j \delta \rho_{Iij} - \mathbb{C}_{JIil} \, \p_l \delta \lambda_J) - C_{Iij}^H(P,L,x,t) \delta\rho_{Iij} \\
    & \qquad \qquad \qquad \qquad - \mathfrak{f}_I(U_H(P,L,x,t)) \delta \lambda_I + \mathfrak{A}_{Ik}(U_H(P,L,x,t)) \p_k \delta \lambda_{I},
    \end{aligned}
\end{equation*}
and collecting terms,
\begin{equation*}
    \begin{aligned} 
    \delta S & = \int_\Omega d^3x \, \delta \lambda_{I}(x, 0) u_I^H(P(x, 0), L(x,0),x , 0) - \int_\Omega d^3x \, \delta \lambda_{I}(x, T) u_I^H(P(x, T), L(x,T),x, T) \\
    & \quad + \int_{\Omega \times (0,T)} d^3xdt  \  \delta \lambda_I \left( \p_t u^H_I(P,L,x,t) - \mathfrak{f}_I(U_H(P,L,x,t)) - \p_k \left( \mathfrak{A}_{Ik}(U_H(P,L,x,t))\right)  \right) \\
    & \quad + \int_{\Omega \times (0,T)} d^3xdt  \  \delta \lambda_I \left( - \mathbb{A}_{IJ} u^H_J (P,L,x,t) - \mathbb{B}_{IJk} \p_k \left(u^H_J(P,L,x,t)\right) - \mathbb{C}_{IJil} \p_l\left( B^H_{Ji}(P,L,x,t) \right) \right) \\
    & \quad + \int_0^T dt \int_{\p \Omega} da \ \delta \lambda_I \left( \mathbb{M}_{jk} \mathfrak{A}_{Ij}(U_H(P,L,x,t)) + u^H_J(P,L,x,t) \mathbb{B}_{IJk} + \mathbb{C}_{IJik} B^H_{JI}(P,L,x,t) \right)\, n_k \\
    & \quad -\int_0^T dt \int_{\p \Omega} da \  \delta \gamma_{Ii} \,u_I^H(P,L,x,t) \,n _i \\
    & \quad + \int_{\Omega \times (0,T)} d^3x dt \, \delta \gamma_{Ii} \left( \p_i u_I^H(P,L,x,t) - B_{Ii}^H(P,L,x,t) \right) \\
    & \quad -\int_0^T dt \int_{\p \Omega} da \  \delta \rho_{Iij} B^H_{Ii}(P,L,x,t) \,  n_j \\
    & \qquad + \int_{\Omega \times (0,T)} d^3x dt \, \delta \rho_{Iij} \left( \p_j B^H_{Ii}(P,L,x,t) - C^H_{Iij}(P,L,x,t)\right).
    \end{aligned}
\end{equation*}
From the above calculation we read off the modification required to $S$ to account for the specified initial and boundary conditions:
\begin{equation}\label{eq:S_ibvp}
\begin{aligned}
    S_{ibvp}[\lambda, \gamma, \rho] & :=  - \int_{\Omega \times (0,T)} d^3x dt\, M^*(P(x,t),L(x,t),x,t) \\
    & \quad - \int_\Omega d^3x \, \lambda_I(x, 0) \, \overline{u}^{(i)}_I (x) - \int_0^T dt \int_{\p \Omega_\tau} da \,\lambda_I(x,t) \,\overline{\tau}_I(x,t)\\
    & \quad + \int_0^T dt \int_{\p \Omega_u} da \, \gamma_{Ii} (x,t) \, \overline{u}^{(b)}_I (x,t)\, n_i(x) + \int_0^T dt \int_{\p \Omega_{\nabla u}} da \, \rho_{Iij}(x,t) \overline{B}_{Ii}(x,t) \, n_j(x)\\
   \mbox{with}  & \quad \lambda_I(x, T) = \mbox{`arbitrarily' prescribed and } \delta \lambda_I(x,T) = 0, x \in \Omega\\
    & \quad \lambda_I(x, t) = \mbox{`arbitrarily' prescribed and } \delta \lambda_I(x,t) = 0, (x,t) \in \p \Omega \backslash \p \Omega_\tau \times (0, T]\\
    & \quad \gamma_{Ii}(x, t) = \mbox{`arbitrarily' prescribed and } \delta \gamma_{Ii}(x,t) = 0, (x,t) \in \p \Omega \backslash \p \Omega_u \times (0, T]\\
    & \quad \rho_{Iij}(x, t) = \mbox{`arbitrarily' prescribed and } \delta \rho_{Iij}(x,t) = 0, (x,t) \in \p \Omega \backslash \p \Omega_{\nabla u} \times (0, T],
\end{aligned}
\end{equation}
(with prescriptions chosen to avoid discontinuities at `space-time corners' of boundary of the space-time domain $\Omega \times (0,T)$).

We end by noting that in ongoing work these ideas have been used to successfully formulate and compute approximate solutions (with minimal error) of the heat equation and the first-order wave equation in bounded domains, in one space dimension and time. These parabolic and hyperbolic equations are solved by a common methodology based on computing weak solutions to degenerate elliptic boundary value problems in a space-time domain, involving oblique natural boundary conditions. Uniqueness of solutions to the corresponding dual problems with the initial-boundary conditions developed following the above ideas is also shown.

%%{\color{red}%%
As already noted in the Introduction, the primary application envisioned for this variational principle is in enabling a Feynman Path Integral based statistical analysis of dislocation dynamics. A second major application is in designing an efficient and robust numerical scheme for the system \eqref{eq:FDM}, which is a \textit{system} of Hamilton-Jacobi equations. It also has the potential of providing a pathway to the rigorous analysis of the system in the hands of bona-fide experts in PDE and variational calculus.%%}

%%%%%%%%%%%%%%%%%%%%%%%%%%%%%
\bibliographystyle{IEEEtran}\bibliography{disloc.bib}
\end{document}